\documentstyle[11pt]{article}
\def\a{\alpha} 
\def\b{\beta} 
\def\g{\gamma}
\def\G{\Gamma}
\def\D{\triangle}
\def\t{\tau}
\def\d{\delta}  
\def\th{\theta} 
\def\l{\lambda}

\def\r{\rho}
\def\ph{\phi}
\def\ps{\psi}
\def\m{\mu}
\def\s{\sigma}
\def\sou{\overline}
\def\so{\underline} 

\def\f{\rightarrow}
\def\fe{\hookrightarrow}
\def\q{\forall}
  
\def\v{\vdash}
 
\def\p{\succ}

\def\R{\ifmmode{\rm I\mkern-3.1mu
R\mkern1mu}\else{\rm I\kern-.18em 
R\hskip1pt\ }\fi\relax} 

\def\Z{\ifmmode{ Z\mkern-4.6mu
Z\mkern2mu}\else{ Z\kern-.28em 
Z\hskip1pt\ }\fi\relax} 

\def\Q{\ifmmode{\rm Q\mkern-10mu
l\mkern4.5mu}\else{\rm Q\kern-.57em
l\hskip3pt\ }\fi\relax} 

\def\N{\ifmmode{\rm I\mkern-3.1mu
N\mkern0.5mu}\else{\rm I\kern-.16em
N\hskip0.5pt\ }\fi\relax} 

\def\C{\ifmmode{\rm C\mkern-8.8mu
l\mkern4mu}\else{\rm C\kern-.48em
l\hskip2.6pt\ }\fi\relax} 

\newtheorem{theo}{Th\'eor\`eme}[section]
\newtheorem{lemma}{Lemme}[section]

\newtheorem{corollary}{Corollaire}[section]

\begin{document} 

\vspace*{0cm}

\Large\bf 
La Valeur d'un Entier Classique en $\l\m$-Calcul \normalsize \rm \\

\bf 
Karim Nour \footnote{ Je remercie M. Parigot pour son aide qui m'a permis de r\'ealiser ce
travail.}\\
\rm
Laboratoire de math\'ematiques, \'equipe de logique,\\ 
Universit\'e de Chamb\'ery, 73376 Le Bourget-du-Lac Cedex, France \footnote{e-mail
nour@univ-savoie.fr} \\

\bf Abstract. \rm
In this paper, we present three methods to give the value of a classical integer in
$\l\m$-calculus. The first method is an external method and gives the value and the false part of a
normal classical integer. The second method uses a new reduction rule and gives as result the
corresponding Church integer. The third method is the M. Parigot's method which uses the J.L.
Krivine's storage operators.\\



\section{Introduction}

Consid\'erons le $\l$-calcul muni d'un syst\`eme de typage bas\'e sur la logique
intuitionniste du second ordre : le syst\`eme $AF2$ de J.L. Krivine. Ce syst\`eme est une simple
extension du syst\`eme $F$ de J.Y. Girard, capable d'exprimer les
sp\'ecifications exactes des programmes, ce qui permet d'obtenir un programme calculant une fonction
en \'ecrivant une d\'emonstration de sa totalit\'e. \\

Comment programme-t-on les fonctions sur les entiers naturels en $AF2$ ?
\begin{itemize}
\item Il faut d'abord exprimer par une formule le type des entiers
naturels : pour ceci on introduit la formule $N[x] = \q X \{ X(0) , \q y(X(y) \f X(sy)) \f X(x) \}$.
Cette formule signifie que $x$ est un entier si et seulement si $x$ appartient au plus petit ensemble
contenant z\'ero et stable par le successeur. La repr\'esentation de l'entier $n$ en $\l$-calcul est
obtenue en d\'emontrant l'\'enonc\'e $N[s^n(0)]$. L'\'el\'ement cl\'e qui fait fonctionner les
choses est l'unicit\'e de la repr\'esentation des entiers. En effet, on d\'emontre facilement que
les entiers de Church sont les seuls $\l$-termes normaux clos de type $N[s^n(0)]$.   
\item Il faut ensuite exprimer les
sp\'ecifications du programme par un syst\`eme d'\'equations d\'efinissant la fonction \`a
calculer.    
\item Il faut finalement d\'emontrer le th\'eor\`eme \'enon\c{c}ant la totalit\'e de la
fonction : si $f$ est une fonction d\'efinie de $\N^r$ dans $\N$, on obtient un programme pour $f$ en
d\'emontrant la formule $\q x_1...\q x_r \{ N[x_1],...,N[x_r] \f N[f(x_1,...,x_r)] \}$. 
\end{itemize}
Pour capturer le contenu algorithmique des preuves classiques, M. Parigot a introduit le
$\l\m$-calcul \'equip\'e d'un syst\`eme de typage bas\'e sur la logique classique du second ordre. Il
a d\'emontr\'e que ce calcul poss\`ede de tr\`es bonnes propri\'et\'es : propri\'et\'e
de Church-Rosser - conservation de type - normalisation forte - ...
Mais la m\'ethode d\'etaill\'ee ci dessus pour programmer des fonctions ne marche pas bien dans ce
syst\`eme. Ceci provient du fait qu'en logique classique on perd la propri\'et\'e de
l'unicit\'e de la repr\'esentation des entiers. De plus il est difficile a priori de conna\^{\i}tre
la valeur d'un entier classique.\\

Nous posons les questions suivantes :
\begin{itemize}
\item Peut-on caract\'eriser les entiers classiques ?
\item Peut-on conna\^{\i}tre la valeur d'un entier classique ?
\end{itemize}

M. Parigot a donn\'e des r\'eponses positives \`a ces questions. En effet, il a trouv\'e un
algorithme qui teste si un $\l \m$-terme normal est un entier classique ou pas, et dans le cas positif
il trouve sa valeur.\\

Si nous voulons construire des programmes par des preuves, nous n'avons pas besoin de v\'erifier si
un $\l \m$-terme normal est un entier ou pas. En effet ceci est assur\'e par le typage. Ce qui nous
int\'eresse le plus, c'est de trouver la valeur d'un entier classique.\\

Dans ce papier, nous pr\'esentons trois m\'ethodes pour trouver la valeur d'un entier
classique. 
\begin{itemize}
\item La premi\`ere m\'ethode est une m\'ethode externe au $\l\m$-calcul et donne, en plus de la
valeur, la partie fausse d'un entier classique normal. En utilisant cette m\'ethode, on peut associer
\`a chaque $\l\m$-terme clos de type $N$ (le type des entiers du syst\`eme de typage $F$) une valeur
fictive.   
\item La deuxi\`eme m\'ethode utilise une nouvelle r\`egle de r\'eduction (la r\`egle de
nettoyage). En ajoutant cette r\`egle au $\l \m$-calcul on obtient un calcul qui poss\`ede des
mauvaises propri\'et\'es (on perd la propri\'et\'e de Church-Rosser et la conservation de type).
Cette m\'ethode donne comme r\'esultat l'entier de Church correspondant. 
\item La troisi\`eme m\'ethode est celle de M. Parigot qui utilise les op\'erateurs de mise en
m\'emoire de J.L. Krivine. Nous pr\'esentons ici des op\'erateurs de mise en m\'emoire qui donnent la
valeur fictive d'un $\l \m$-terme de type $N$ et nous montrons que ce n'est pas toujours le cas
pour les autres op\'erateurs.  
\end{itemize}

\section{$\l\m$-calcul}

Dans ce paragraphe, nous pr\'esentons la d\'eduction naturelle classique du second ordre ainsi que
son interpr\'etation calculatoire, le $\l\m$-calcul. Nous utilisons un syst\`eme de
d\'eduction naturelle avec plusieurs conclusions. Le $\l\m$-calcul est une simple extension du
$\l$-calcul qui donne exactement le contenu algorithmique des preuves \'ecrites dans ce syst\`eme.

\subsection{Le $\l\m$-calcul pur}

Le $\l\m$-calcul poss\`ede deux alphabets distincts de variables : un ensemble de $\l$-variables
$x,y,z$,..., et un ensemble de $\m$-variables $\a,\b,\g$,....  \\ 
Les termes sont d\'efinis de la mani\`ere inductive suivante : 
\begin{itemize}
\item[] - Si $x$ est une $\l$-variable, alors $x$ est un terme ;
\item[] - Si $x$ est une $\l$-variable et $u$ est un terme, alors $\l x u$ est un terme ;
\item[] - Si $u$ et $v$ sont des termes, alors $(u)v$ est un terme; 
\item[] - Si $t$ est un terme et $\a,\b$ sont des $\m$-variables, alors $\m\a[\b]t$ est un terme.
\end{itemize}
Les termes du $\l$-calcul sont obtenus seulement par les trois premi\`eres clauses. \\ 
Un terme nomm\'e est une expression de la forme $[\a]t$ o\`u $t$ est un terme et $\a$ est une
$\m$-variable. Dans ce cas on dit que le terme $t$ est nomm\'e par $\a$. On consid\`ere parfois 
les termes nomm\'es comme des termes. Pour simplifier, on suppose que l'ensemble des variables
libres et l'ensemble des variables li\'ees d'un terme sont toujours distincts, et qu'une
variable li\'ee est li\'ee une seule fois.\\
Les termes du $\l\m$-calcul sont appel\'es $\l\m$-termes, et les termes du $\l$-calcul,
$\l$-termes.\\ 

La r\'eduction en $\l\m$-calcul est induite par deux notions diff\'erentes de r\'eductions :\\ 

\bf Les r\`egles de calcul : \/ \rm
\begin{itemize} 
\item[] ($C_1$) $(\l xu)v \f u[v/x]$ (la $\b$-r\'eduction)
\item[] ($C_2$) $(\m\a u)v \f \m\a u[v/$*$\a]$ (la $\m$-r\'eduction)
\item[] o\`u $u[v/$*$\a]$ est obtenu \`a partir du $u$ en rempla\c{c}ant inductivement chaque
sous terme de la forme $[\a]w$ par $[\a](w)v$.  
\end{itemize} 
\bf Les r\`egles de simplification : \/ \rm
\begin{itemize} 
\item[] ($S_1$)  $[\a]\m\b u \f  u[\a/\b]$
\item[] ($S_2$)  $\m\a[\a]u \f  u$, si $\a$ n'a pas d'occurrences libres dans $u$   
\item[] ($S_3$) $\m\a u \f \l x\m\a u[x/$*$\a]$, si $u$ contient un sous terme de la
forme $[\a]\l yw$.   \end{itemize}
Il est clair que les r\`egles de simplifications sont fortement normalisables, en effet, les
r\`egles $(S_i)_{1\leq i\leq 3}$ ne cr\'eent pas de nouveaux r\'edex ou diminuent strictement la
longueur d'un $\l\m$-terme.

\begin{theo}[Th\'eor\`eme de Church-Rosser]  
En $\l\m$-calcul, la r\'eduction est confluente (c'est \`a dire si $u \f u_1$ et $u \f
u_2$, alors il existe $v$ tel que $u_1 \f v$ et $u_2 \f v$).
\end{theo}
\bf Preuve \rm Voir [4]. $\Box$ \\ 

L'ensemble des $\l\m$-termes en formes normales de t\^ete $H$ est d\'efini de la
mani\`ere inductive suivante :
\begin{itemize}
\item[] - Si $u \in H$, alors $\l xu \in H$ ;
\item[] - Si $u \in H$, et $(u)v$ n'est pas un r\'edex, alors $(u)v \in H$ ;
\item[] - Si $u \in H$, et $\m\a[\b]u$ n'est pas un r\'edex, alors $\m\a[\b]u \in H$.
\end{itemize} 

Une r\'eduction de t\^ete, est une suite de r\'eductions $t_1$,...,$t_n$ telle que pour tout
$1\leq i\leq n-1$, $t_{i+1}$ est obtenu en r\'eduisant le r\'edex le plus \`a gauche de
$t_i$ et $t_i$ n'est pas en forme normale de t\^ete. Si $u$ est obtenu \`a partir de
$t$ par une r\'eduction de t\^ete, on note $t \p u$. \\
Si $t\p u$, on note $h(t,u)$ la longueur de la r\'eduction de t\^ete entre $t$ et $u$. \\
L'\'equivalence de t\^ete est not\'ee : $u\sim v$ ssi il existe $w$ tel que $u\p w$ et $v\p w$. \\

Il est facile de prouver le lemme suivant (ce lemme peut \^etre prouv\'e par induction, comme dans le
cas du $\l$-calcul (voir [2])).

\begin{lemma} Si $u\p v$, alors :\\
1) $u[\sou{p} / \sou{x} , \sou{q} /$*$ \sou{\a}]\p v[\sou{p} / \sou{x} , \sou{q} /$*$ \sou{\a}]$
et $h(u[\sou{p} / \sou{x} , \sou{q}  /$*$ \sou{\a}],v[\sou{p} / \sou{x} , \sou{q} /$*$
\sou{\a}])=h(u,v)$. \\ 
2) Pour tout $\sou{w}$, il existe $w$ tel que $(u)\sou{w} \p w$, $(v)\sou{w}\p w$, et
$h((u)\sou{w},w)=h((v)\sou{w},w)+h(u,v)$. 
\end{lemma}

\bf Remarque. \rm Le lemme 2.1 montre que pour effectuer la r\'eduction de t\^ete de
$u[\sou{p} / \sou{x} , \sou{q} $/*$ \sou{\a}]$ (resp. $(u)\sou{w}$), il est \'equivalent (m\^eme r\'esultat, et
m\^eme nombre de pas) d'effectuer un certain nombre de pas dans la r\'eduction de t\^ete de $u$ pour
obtenir $v$, puis de faire la r\'eduction de t\^ete de $v[\sou{p} / \sou{x} , \sou{q} $/*$
\sou{\a}]$ (resp. $(v)\sou{w}$).

\subsection{Le $\l\m$-calcul typ\'e}

Les types sont les formules du calcul des pr\'edicats du second ordre.\\ 
Les connecteurs logiques utilis\'es sont $\perp$, $\f$ et $\q$.\\
Les langages contiennent des variables d'individus (ou du premier ordre) not\'ees $x,y,z$,..., et
des variables de relations (ou du second ordre) not\'ees $X,Y,Z$,....\\
Nous ne supposons pas que le langage contient un symbole de constante sp\'ecial pour l'\'egalit\'e.
Cependant, on d\'efinit la formule $u=v$ (o\`u $u,v$ sont des termes) par $\q Y(Y(u)\f Y(v))$, o\`u
$Y$ est une variable de relation unaire. Une telle formule est appel\'ee \'equation. On note
$\v_E u=v$, si $u=v$ est une cons\'equence de l'ensemble d'\'equations $E$. \\  
La formule $F_1 \f (F_2 \f (...\f (F_n \f G)...))$ est not\'ee
$F_1,F_2,...,F_n \f G$ et la formule $F \f \perp$ est not\'ee $\neg F$.\\  
Les preuves sont \'ecrites dans un syst\`eme de d\'eduction naturelle avec plusieurs conclusions
pr\'esent\'ees avec des s\'equents : 
\begin{itemize}
\item[] - Les formules \`a gauche de $\v$ sont \'etiquet\'ees par des $\l$-variables ;  
\item[] - Les formules \`a droite de $\v$ sont \'etiquet\'ees par des $\m$-variables, except\'e
d'une formule qui est \'etiquet\'ee par un $\l\m$-terme ;  
\item[] - Des formules distinctes ne poss\`edent pas la m\^eme \'etiquette.
\end{itemize}
Les parties gauches et droites d'un s\'equent sont consid\'er\'ees comme des ensembles,
et donc la contraction des formules est donn\'ee implicitement. Les affaiblissements sont inclus
dans les r\`egles (2) et (9).\\
Soient $t$ un $\l\m$-terme, $A$ un type, $\G=x_1 :A_1,...,x_n :A_n$,
et $\D=\a_1 :B_1,...,\a_m :B_m$ deux contextes. On d\'efinit par les r\`egles suivantes
la notion "t est de type A dans $\G$ et $\D$", et on \'ecrit $\G\v t:A,\D$. 
\begin{itemize} 
\item[] (1) $\G\v x_i :A_i,\D$ $1\leq i\leq n$
\item[] (2) Si $\G,x:A \v t:B,\D$, alors $\G\v \l xt:A\f B,\D$  
\item[] (3) Si $\G\v u:A\f B,\D$, et $\G\v v:A,\D$, alors $\G\v (u)v:B,\D$
\item[] (4) Si $\G\v t:A,\D$, alors $\G\v t:\q xA,\D$ (*) 
\item[] (5) Si $\G\v t:\q xA,\D$, alors $\G\v t:A[u/x],\D$ (**)   
\item[] (6) Si $\G\v t:A,\D$, alors $\G\v t:\q XA,\D$ (*) 
\item[] (7) Si $\G\v t:\q XA,\D$, alors $\G\v t:A[G/X],\D$ (**) 
\item[] (8) Si $\G\v t:A[u/x],\D$, alors $\G\v t:A[v/x],\D$ (***) 
\item[] (9) Si $\G\v t:A,\b:B,\D$, alors :
\begin{itemize} 
\item[] -- $\G\v \m\b[\a]t:B,\a :A,\D$ si $\a \neq \b$
\item[] -- $\G\v \m\a[\a]t:B,\D$ si $\a = \b$
\end{itemize} 
\end{itemize} 
Les r\`egles pr\'ec\'edentes font l'object des restrictions suivantes :\\
(*) Les variables $x$, et $X$ n'ont pas d'occurences libres dans $\G$.\\
(**) $u$ est un terme et $G$ est une formule du langage.\\
(***) $u$ et $v$ sont des termes, tels que $u=v$ est une cons\'equence d'un ensemble d'\'equations.
 
\begin{theo}[Th\'eor\`emes de conservation de type et de normalisation forte]. \\
1) Le type est pr\'eserv\'e durant une r\'eduction (c'est \`a dire si $\G\v u:A,\D$, et
$u\f v$, alors $\G\v v:A,\D$).\\
2) Les $\l\m$-termes typables sont fortement normalisables (c'est \`a dire si $\G\v
t:A,\D$, alors toute r\'eduction qui commence par $t$ est fini).
\end{theo}
\bf Preuve \rm Voir [4] et [6]. $\Box$  

\section{Les entiers}

\subsection{Les entiers intuitionnistes}

Si on oublie la r\`egle (9), on obtient le syst\`eme de J.L. Krivine appel\'e $AF2$ (voir [1]). Dans
ce syst\`eme on utilise les $\l$-termes et on garde uniquement la notion de la $\b$-r\'eduction.\\
Dans le syst\`eme de typage $AF2$, chaque type de donn\'ees peut \^etre d\'efini par une formule. Par
exemple le type des entiers est la formule :
$N[x]=\q X\{X(0) , \q y(X(y) \f X(sy)) \f X(x) \}$ o\`u $X$ est une variable de relation unaire, $0$
est un symbole de constante pour le z\'ero, et $s$ est un symbole de fonction unaire pour le
successeur.\\ 
Le $\l$-terme $\so{0}=\l x\l fx$ est de type $N[0]$ et repr\'esente "z\'ero".\\
Le $\l$-terme $\so{s}=\l n\l x\l f(f)((n)x)f$ est de type $\q y(N[y] \f N[s(y)])$ et
repr\'esente la fonction "successeur".\\  
Un ensemble d'\'equations $E$ est dit ad\'equat pour le type des entiers ssi :
\begin{itemize} 
\item[] - $\not\v_E s(a)=0$ ;
\item[] - Si $\v_E s(a)=s(b)$, alors $\v_E a=b$.
\end{itemize}
Dans la suite, on suppose que tous les ensembles d'\'equations sont ad\'equats pour le type
des entiers.

\begin{theo}[Unicit\'e de la repr\'esentation des entiers]                         
Pour tout entier $n$, $\so{n}=\l x\l f(f)^nx$ est l'unique $\l$-terme normal clos de type $N[s^n(0)]$.
\end{theo}
\bf Preuve \rm Voir [1]. $\Box$ \\ 

La trace propositionnelle $N=\q X\{X , X \f X \f X \}$ de $N[x]$ d\'efinit aussi les entiers.

\begin{theo} Un $\l$-terme normal clos est de type $N$ ssi il est de la forme $\so{n}$, pour un
certain entier $n$.
\end{theo}
\bf Preuve \rm Voir [1]. $\Box$ \\  

\begin{theo}[Th\'eor\`eme de programmation]  
Soit $F$ un $\l$-terme clos de type 
\begin{displaymath}
\q x_1...\q x_r (N[x_1],...,N[x_r] \f N[f(x_1,...,x_r)])
\end{displaymath}
modulo un ensemble d'\'equations $E$ ad\'equat pour le type des entiers. Alors $F$ est un programme
pour $f$ dans le sens suivant : pour tout $n_1,...,n_r,m$, $(F)\so{n_1}...\so{n_r} \f \so{m}$ ssi
$\v_E f(s^{n_1}(0),...,s^{n_r}(0))=s^m(0)$.
\end{theo}
\bf Preuve \rm Voir [1]. $\Box$ 

\subsection{Les entiers classiques}

En d\'eduction naturelle classique l'unicit\'e de la repr\'esentation des entiers ne reste pas vraie.
En effet l'affaiblissement \`a droite cr\'ee des (fausses) copies des entiers, et la
contraction droite permet de regrouper des copies diff\'erentes dans une seule.\\

Un entier classique est un $\l\m$-terme clos $\th$, tel que $\v \th :N[s^n(0)]$ pour un certain
entier $n$. On dit aussi que l'entier classique $\th$ est de valeur $n$.\\

On va carat\'eriser maintenant les entiers classiques.

\section{L'algorithme de Parigot caract\'erisant les entiers classiques }

Soient $x$ et $f$ deux variables fixes, et $N_{x,f}$ l'ensemble des $\l\m$-termes d\'efinis par :

\begin{itemize}
\item[] - $x\in N_{x,f}$ ;
\item[] -  $\m\a[\b]x\in N_{x,f}$ ;
\item[] - Si $u\in N_{x,f}$, alors $(f)u\in N_{x,f}$ ;
\item[] - Si $u\in N_{x,f}$, alors $\m\a[\b](f)u\in N_{x,f}$.
\end{itemize}

\begin{theo}
Les $\l\m$-termes normaux clos de type $N$ sont exactement les $\l\m$-termes de la forme $\l x\l fu$
o\`u $u\in N_{x,f}$ sans $\m$-variables libres. 
\end{theo} 
\bf Preuve \rm Voir [5]. $\Box$ \\

On d\'efinit, pour chaque $u\in N\sb{x,f}$, l'ensemble $rep(u)$, qui d\'esigne intuitivement
l'ensemble des entiers repr\'esent\'es par $u$ :
\begin{itemize}
\item[] - $rep(x)=\{ 0 \}$ ;
\item[] - $rep((f)u)=\{ n+1$ si $n\in rep(u) \}$ ;
\item[] - $rep(\m\a[\b]u)=\bigcap rep(v)$ pour chaque sous terme $[\a]v$ de $[\b]u$.
\end{itemize}

Le th\'eor\`eme suivant caract\'erise les entiers classiques.

\begin{theo} Les entiers classiques normaux de valeur $n$ sont exactement les
$\l\m$-termes de la forme $\l x\l fu$ o\`u $u\in N_{x,f}$ sans $\m$-variables libres et
tel que $rep(u)=\{ n \}$. \end{theo}
\bf Preuve \rm Voir [5]. $\Box$ \\

\bf Exemples \\ \rm
1) Il est facile de v\'erifier que $\so{0}$ est l'unique entier classique normal de valeur $0$.\\
2) Il est \'egalement facile de v\'erifier que tout entier classique de valeur $1$ est
\'equivalent \`a  un $\l\m$-terme de la forme : 
\begin{displaymath}
\l x \l f \m \a [\a](f)\m \a_1[\a](f)...\m\a_n[\a](f)\m \b[\a](f)\m \b_1[\a](f)...\m \b_{m-1}[\a](f)
\m \b_m[\b]x 
\end{displaymath}
o\`u $n,m \geq 0$, et $\a,\b,\a_i,\b_j$ avec $1\leq i\leq n$ et $1\leq j\leq m$ sont des
$\m$-variables distinctes.\\
3) Soit $\th=\l x\l
f(f)\m\a[\a](f)\m\ph[\a](f)\m\ps[\a](f)(f)\m\b[\ph](f)\m\d[\b](f)\m\g[\a](f)\m\r[\b](f)x$.\\
$\th$ est-il un entier classique ? si oui quelle est sa valeur ?\\
Posons $u=(f)\m\a[\a](f)\m\ph[\a](f)\m\ps[\a](f)(f)\m\b[\ph](f)\m\d[\b](f)\m\g[\a](f)\m\r[\b](f)x$.\\
On v\'erifie que $rep(u)=\{ 4 \}$, donc $\th$ est un entier classique de valeur 4.\\
Remarquons que si nous savons \`a l'avance que $\th$ est un entier classique, alors nous n'avons
pas besoin de chercher $rep(u)$ pour tous les sous termes de $u$. Dans ce cas, il suffit de s'occuper
des sous termes de $\th$ qui repr\'esentent un seul entier.

\section{Une m\'ethode simple pour trouver la valeur d'un entier classique}

On pr\'esente maintenant une m\'ethode simple pour trouver la valeur d'un entier classique.\\

On d\'efinit, pour chaque $u\in N\sb{x,f}$, l'ensemble $val(u)$, qui d\'esigne intuitivement
l'ensemble des valeurs possibles de $u$ : 
\begin{itemize} 
\item[] - $val(x)=\{ 0 \}$ ; 
\item[] - $val((f)u)=\{ n+1$ si $n\in val(u) \}$ ;
\item[] - $val(\m\a[\b]u)=\bigcup val(v)$ pour chaque sous terme $[\a]v$ de $[\b]u$.
\end{itemize}

Soient $u\in N\sb{x,f}$ sans $\m$-variables libres et $\a_1,...,\a_n$ les $\m$-variables de $u$ qui
v\'erifient :
\begin{itemize}  
\item[] - $\a_1$ est la $\m$-variable, telle que $[\a_1](f)^{i_1 }x$ est un sous terme de $u$ 
\item[] - $\a_j$ $2\leq j\leq n$ est la $\m$-variable, telle que $[\a_j](f)^{i_j}\m\a_{j-1}u_{j-1}$
est un sous terme de $u$   
\item[] - $u=(f)^{i_{n+1}}\m\a_n u_n$.  
\end{itemize}

Sch\'ematisons cette derni\`ere d\'efinition par le dessin suivant :\\

$ u = (f)^{i_{n+1}}\m\a_n 
\underbrace{...[\a_3](f)^{i_3}\m\a_2
\underbrace{...[\a_2](f)^{i_2}\m\a_1
\underbrace{...[\a_1](f)^{i_1}x}_{u_1}
                                      }_{u_2}
                                             }_{u_n} $.\\

Soient $t_0 =x$ et $t_j=\m\a_j u_j$ $1\leq j \leq n$. 

\begin{lemma}
Pour tout $1\leq j \leq n+1$ on a :\\  
1) $val(t_{j-1})=\displaystyle{\{\sum_{1 \leq k \leq j} \: i_k \}}$.\\ 
2) Pour chaque sous terme $t$ de $u_j$, tel que $t \neq (f) ^r t_k$  $0\leq k\leq j-1$, on a
$val(t)=\emptyset$. \\
En particulier $val(u)= \displaystyle{\{ \sum_{1 \leq k \leq n+1} \: i_k \}} $.
\end{lemma}
\bf Preuve \rm Par induction sur $j$.
\begin{itemize}  
\item[] - Pour $j=1$ : on a $val(t_0)=\{ 0 \}$.\\
Soit $t$ un sous terme de $u_1$, tel que $t\neq (f)^r t_0$. D\'emontrons, par induction sur $t$, que
$val(t)=\emptyset$.
\begin{itemize}  
\item[] - Si $t=(f)t'$, alors, par hypoth\`ese d'induction, $val(t')=\emptyset$, et
$val(t)=\emptyset$.
\item[] - Si $t=\m\b t'$, alors, il est clair que $\b\neq \a_1$. Donc, par hypoth\`ese
d'induction, pour tout terme $u$ nomm\'e $\b$, on a $val(u)=\emptyset$. Donc $val(t)=\emptyset$.
\end{itemize}  
\item[] - Supposons la propri\'et\'e vraie pour tout $1\leq i\leq j-1$, et montrons la
pour $j$.\\ 
Par d\'efinition, on a $val(t_j)=\bigcup val(t)$ pour chaque sous terme $[\a_j]t$ de $u_j$. Donc, par
hypoth\`ese d'induction, $val(t_j)=\displaystyle{\{ i_j + \sum_{1 \leq k \leq j-1} \: i_k \}}
\bigcup \emptyset = \displaystyle{\{ \sum_{1 \leq k \leq j} \: i_k \}}$.\\
Soit $t$ un sous terme de $u_j$, tel que $t\neq (f)^r t_k$  $0\leq k\leq j-1$. D\'emontrons, par
induction sur $t$, que $val(t)=\emptyset$.
\begin{itemize}  
\item[] - Si $t=(f)t'$, alors, par hypoth\`ese d'induction, $val(t')=\emptyset$, et
$val(t)=\emptyset$.
\item[] - Si $t=\m\b t'$, alors, il est clair que $\b\neq \a_i$ $1\leq i\leq j-1$. Donc, par
hypoth\`ese d'induction, pour tout terme $u$ nomm\'e $\b$, on a $val(u)=\emptyset$. Donc
$val(t)=\emptyset$. $\Box$ \end{itemize} 
\end{itemize}

\bf Remarque \rm D'apr\`es le lemme 5.1,  on peut associer \`a chaque $\l\m$-terme normal clos
$\th = \l x\l fu$ de type $N$ l'entier $val(u)$. Cet entier repr\'esente la valeur "fictive" de
$\th$.

\begin{lemma} 
Pour chaque $u\in N_{x,f}$, on a $rep(u) \subseteq val(u)$.
\end{lemma}
\bf Preuve \rm Trivial, par d\'efinition de $rep(u)$ et $val(u)$. $\Box$ 

\begin{theo}
Si $\th$ est un entier classique normal de valeur $n$, alors $\th=\l x\l fu$ o\`u $u\in N_{x,f}$ sans
$\m$-variables libres et tel que $val(u)=\{ n \}$. \end{theo}
\bf Preuve \rm D'apr\`es le th\'eor\`eme 4.2, $\th=\l x\l fu$ o\`u $u\in N\sb{x,f}$
sans $\m$-variables libres et tel que $rep(u)=\{ n \}$. D'apr\`es le lemme 5.1, on a
$val(u)=\{ m \}$, et d'apr\`es le lemme 5.2, $rep(u)=\{ n \} \subseteq \{m\} =val(u)$.
Donc $val(u)=\{ n \}$. $\Box$\\

\bf Conclusion \rm Pour trouver la valeur d'un entier classique normal $\th = \l x\l f u$, on parcourt
le $\l\m$-terme $u$ du droite \`a gauche pour chercher les $\m$-variables $\a_j$ $1\leq j\leq n$ et
les entiers $i_j$ $1\leq j\leq n+1$ d\'efinis avant le lemme 5.1. D'apr\`es le lemme 5.1 et le
th\'eor\`eme 5.1, la valeur de l'entier $\th$ est \'egale \`a $\displaystyle{ \sum_{1 \leq k \leq
n+1} \: i_k } $.\\ 

La vraie partie de l'entier de $\th$ est le $\l$-terme $\l x\l
f(f)^{i_{n+1}}(f)^{i_n}...(f)^{i_1}x$.\\

\bf Exemple \rm Reprenons le dernier exemple du paragraphe 4.\\

$ \th = \l x\l f \overbrace{
(f)}^{4}\underbrace{\m\a[\a](f)\m\ph[\a](f)\m\ps[\a]}_{partie\quad fausse}\overbrace{(f)(f)}^{3}
\underbrace{\m\b[\ph](f)\m\d[\b](f)\m\g[\a](f)\m\r[\b]}_{partie\quad fausse}\overbrace{(f)x}^{1} $.\\

On v\'erifie facilement que :
\begin{itemize}
\item[] - $\a_1 = \b$, et $\a_2 =\a$ ;
\item[] - $i_1 = 1$, $i_2 = 2$, et $i_3 = 1$.
\end{itemize}
Donc la valeur de $\th$ est \'egale \`a 4.

\section{Une nouvelle r\`egle de r\'eduction pour trouver la valeur d'un entier classique} 

Nous ajoutons au $\l\m$-calcul une nouvelle r\`egle de r\'eduction (r\`egle de nettoyage) :
\begin{itemize}
\item[] ($N$) $(u)\m\a v \f \m\a v[u/$**$\a]$  
\item[] o\`u $v[u/$**$\a]$ est obtenu \`a partir du $v$ en repla\c{c}ant inductivement chaque
sous terme de la forme $[\a]w$ par $[\a](u)w$. 
\end{itemize}
Nous appelons $\l\m'$-calcul ce nouveau calcul.\\

Le $\l\m'$-calcul poss\`ede des mauvaises propri\'et\'es :
\begin{itemize} 
\item La r\'eduction n'est pas confluente.
\item Le type n'est pas pr\'eserv\'e durant une r\'eduction.
\end{itemize}

Par contre on ne sait pas, pour le moment, si en $\l\m'$-calcul le th\'eor\`eme de normalisation
forte reste valable.\\

On note $u \fe v$, si $u$ se r\'eduit \`a $v$ en utilisant, \`a chaque \'etape la r\`egle de nettoyage
suivie par toutes les r\`egles de simplifications possibles. Cette r\'eduction est bien d\'efinie car
les r\`egles de simplifications sont fortement normalisables.\\

Un $\l\m$-terme est dit normal s'il est normal en $\l\m$-calcul (c.\`a.d. sans tenir compte de la
r\`egle de nettoyage).\\

Nous pr\'esentons maintenant comment utiliser la r\`egle $(N)$ pour trouver la valeur d'un
entier classique.\\

Chaque $u \in N_{x,f}$ est de la forme 
\begin{displaymath}
u = (f)^{j_{m+1}} \m \b_m [\g_m] (f)^{j_m}...\m
\b_2 [\g_2](f)^{j_2} \m \b_1 [\g_1] (f)^{j_1} x 
\end{displaymath}
On associe \`a $u$ la liste des entiers $L(u)=(j_{m+1},...,j_1)$, l'entier $l(u)=m$, et le couple
$C(u)=(l(u),L(u))$. On ordonne $L(u)$ et $C(u)$ lexicographiquement.

\begin{theo}
Soit $u \in N_{x,f}$ sans $\m$-variables libres, tel que $l(u) \not = 0$. Alors il existe
$v \in N_{x,f}$ sans $\m$-variables libres tel que $u \fe v$, $C(u) > C(v)$ et $val(u)=val(v)$.
\end{theo}

\bf Preuve \rm Comme $l(u) \not = 0$, alors $u$ contient un sous terme $u'=(f)^{j-1}(f)
\m \a [\b] w$.\\ 
$u' \fe v'=(f)^{j-1} \m \a [\b] w [f /$**$\a]$, donc on distingue deux cas :
\begin{itemize}
\item[] - Si $j > 1$, soit $v$ le $\l\m$-terme obtenu en rempla\c{c}ant dans $u$ le sous
terme $u'$ par $v'$. Il est clair que $l(u)=l(v)$ et $L(u) > L(v)$, donc $C(u) > C(v)$. 
\item[] - Si $j = 1$, on distingue deux sous cas :
\begin{itemize}
\item[] - Si $u = u'$, posons $v = v'$. Il est clair que $l(u)=l(v)$ et $L(u) > L(v)$, donc $C(u) >
C(v)$.  
\item[] - Si $\m \g [\d] u'$ est un sous terme de $u$, soit $v$ le $\l \m$-terme obtenu en
remplacant dans $u$ le sous terme $\m \g [\d] u'$ par $\m \g [\b]v'[ \d / \a]$ si $\g \not = \b$ et
$v'[ \d / \a]$ si $\g = \b$ ne figure pas dans $v'[ \d / \a]$. Il est clair que $u \fe v$ et $l(u)
> l(v)$, donc $C(u) > C(v)$.   \end{itemize}
\end{itemize}
Dans, tous ces cas, on remarque que $v \in N_{x,f}$ est sans $\m$-variables libres. \\
De plus, on a $val(u)=val(v)$, en effet, il faut distinguer deux cas :
\begin{itemize}
\item[] - Si la $\m$-variable $\a$ est l'une des variables $\a_k$ (d\'efinies avant le lemme 5.1),
alors dans $v$ on a : $j=i_{k+1}-1$ et le $-1$ est compons\'e par la pr\'esense de la variable $f$
devant tous les crochets
 $[\a]$. 
\item[] - Si la $\m$-variable $\a$ n'est pas l'une des variables $\a_k$ (d\'efinies avant le lemme
5.1), alors les entiers $i_k$ restent les m\^emes. $\Box$
\end{itemize}

\bf Remarque \rm La preuve du th\'eor\`eme 6.1 montre que la r\'eduction $\fe$ est fortement
normalisable sur les \'el\'ements de l'ensemble $N_{x,f}$.

\begin{corollary}
Si $u \in N_{x,f}$ sans $\m$-variables libres tel que $val(u) = n$, alors $u \fe (f)^n(x)$.
\end{corollary}

\bf Preuve \rm  D'apr\`es le th\'eor\`eme 6.1, $u \fe \l x\l f v$ tel que $l(v)=0$ et
$val(v)=\{ n \}$, donc $v=(f)^n x$. $\Box$
 
\begin{corollary} \footnote{Ce r\'esultat a \'et\'e cit\'e dans un article de M. Parigot (voir [5])}
Si $\th$ est un entier classique normal de valeur $n$, alors $\th \fe \so{n}$.
\end{corollary}

\bf Preuve \rm D'apr\`es le th\'eor\`eme 4.2, chaque entier classique normal $\th$ est de la forme
$\th = \l x \l f u$ o\`u $u \in N_{x,f}$ sans $\m$-variables libres et $val(u)=\{ n \}$. Donc,
d'apr\`es le corollaire 6.1, $\th \fe \so{n}$. $\Box$\\

\bf Conclusion \rm Pour trouver la valeur d'un entier classique normal $\th$, on
r\'eduit $\th$ avec la nouvelle r\`egle de r\'eduction ($\fe$). D'apr\`es le corollaire 6.2, $\th \fe
\so{n}$ o\`u $n$ est la valeur de $\th$.\\ 

\bf Exemple \rm Reprenons le dernier exemple du paragraphe 4.\\
$ \th = \l x\l f
(f)\m\a[\a](f)\m\ph[\a](f)\m\ps[\a](f)(f)\m\b[\ph](f)\m\d[\b](f)\m\g[\a](f)\m\r[\b](f)x $.\\
On v\'erifie facilement que $\th \fe \so{4}$ (effectuer, par exemple, la r\'eduction de la droite vers
la gauche).

\section {Les op\'erateurs de mise en m\'emoire pour trouver la valeur d'un entier classique}

Un $\l$-terme clos $T$ est dit op\'erateur de mise en m\'emoire pour les entiers ssi
pour tout $n\geq 0$, il existe un $\l$-terme $\t_n\simeq\sb{\b}\so{n}$, tel que pour
tout $\th_n\simeq\sb{\b}\so{n}$, il existe une substitution $\s$, telle que
$(T)\th_nf \p (f)\s(\t_n)$.\\

\bf Exemples\/ \rm 
Posons :\\ $T_1=\l n((n)\d)G$ o\`u
$G=\l x\l y(x)\l z(y)(\so{s})z$ et $\d=\l f(f)\so{0}$ ;\\
$T_2=\l n\l f(((n)f)F)\so{0}$
o\`u $F=\l x\l y(x)(\so{s})y$. \\
Il est facile de v\'erifier que : \\
pour tout $\th_n\simeq\sb{\b}\so{n}$,
(($T_1)\th_n)f \p (f)(\so{s})^n\so{0}$ et
(($T_2)\th_n)f \p (f)(\so{s})^n\so{0}$ (voir [2] et [3]). \\
Donc $T_1$ et $T_2$ sont des op\'erateurs de mise en m\'emoire pour les entiers.\\

Soit $N$*$[x]=\q X\{\neg X(0) , \q y(\neg X(y) \f \neg X(sy)) \f \neg X(x) \}$.\\
Il est facile de v\'erifier que $\v_{AF2} T_1,T_2:\q x\{N$*$[x] \f \neg\neg N[x] \}$ (voir [2] et
[3]).

\begin{theo}
Si $\v_{AF2}$T:$\q$x$\{$N*[x]$\f\neg\neg$N[x]$\}$, alors T est un op\'erateur de mise en m\'emoire
pour les entiers.
\end{theo} \bf Preuve \rm Voir [2]. $\Box$

\begin{theo}
Soient $u\in N_{x,f}$ sans $\m$-variables libres, tel que $val(u)=\{n\}$, et $\th=\l x\l
fu$. Alors $((T_1)\th)f \p (f)(\so{s})^{n}\so{0}$. 
\end{theo}  
\bf Preuve \rm La preuve de ce th\'eor\`eme est assez technique. Nous pr\'esentons
seulement les grandes lignes de la d\'emonstration.\\ 
Reprenons les notations du lemme 5.1, et notons
\begin{itemize}
\item[] - $\sou{\m}$ la suite des $\m$-variables de $u$ ; 
\item[] - $\sou{\a}$ la suite des $\m$-variables $\a_1,...\a_n$ ; 
\item[] - $\sou{\b}$ la suite des $\m$-variables $\sou{\m}$- $\sou{\a}$. 
\end{itemize}
Soit $r_j=\displaystyle{\sum_{j \leq k \leq n+1} \: i_k}$ et $S_j=(\so{s})^{r_j}\so{0}$ $1\leq
j\leq n$. Notons $\sou{S}=S_1,...,S_n$.\\
On prouve par induction sur $n-j$ que : \\ 
Pour chaque $1\leq j\leq n$, on a $((T_1)\th)f \sim \m\a[\a]\t_j[\a/\sou{\m}]$ o\`u \\ 
$\t_j=t_j[f/x,F/f,\sou{S}/$*$\sou{\a},\sou{S'}$/*$\sou{\b}](\so{s})^{r_j}\so{0}$ et les
\'el\'ements de $\sou{S'}$ sont de la forme $(\so{s})^{r}\so{0}$.\\  
Donc en particulier, $((T_1)\th)f \sim (f)(\so{s})^{n}\so{0}$. $\Box$ 

\begin{theo}
Soient $u\in N_{x,f}$ sans $\m$-variables libres, tel que $val(u)=\{n\}$, et $\th=\l x\l
fu$. Alors $((T_2)\th)f \p (f)(\so{s})^{n}\so{0}$. 
\end{theo}  
\bf Preuve \rm M\^eme preuve que celle du th\'eor\`eme 7.2. $\Box$ 

\begin{corollary}
Si $\th$ est un entier classique normal de valeur $n$, alors
$((T_i)\th) f \p (f)(\so{s})^{n}\so{0}$ [i=1 ou 2].
\end{corollary}  
\bf Preuve \rm On utilise les th\'eor\`emes 7.2 et 7.3. $\Box$ \\

\bf Conclusion \rm Pour trouver la valeur d'un entier classique normal $\th$, on
r\'eduit $((T_i)\th)\l xx$ [i=1 ou 2]. D'apr\`es le th\'eor\`eme 7.3, $((T_i)\th) \l xx \p
(\l xx)(\so{s})^{n}\so{0} \f \so{n}$ o\`u $n$ est la valeur de $\th$.\\ 

\bf Exemple \rm Reprenons le dernier exemple du paragraphe 4.\\
$ \th = \l x\l f
(f)\m\a[\a](f)\m\ph[\a](f)\m\ps[\a](f)(f)\m\b[\ph](f)\m\d[\b](f)\m\g[\a](f)\m\r[\b](f)x $.\\
On v\'erifie facilement que $((T_i)\th)\l xx \f \so{4}$ [i=1 ou 2]. \\

Le corollaire 7.1 est un cas particulier du th\'eor\`eme suivant :
\begin{theo}
Si $\v_{AF2}$T:$\q$x$\{$N*[x]$\f\neg\neg$N[x]$\}$, alors pour tout $n\geq 0$, il existe un $\l$-terme
$\t_n\simeq\sb{\b}\so{n}$, tel que pour tout entier classique normal $\th$ de valeur $n$, il existe
une substitution $\s$, telle que $((T)\th)f \sim \m\a[\a](f)\s(\t_n)$.  
\end{theo}  
\bf Preuve \rm Voir [5]. $\Box$ \\

\bf Remarque \/ \rm
Le th\'eor\`eme 7.2 ne reste pas vrai pour un op\'erateur de mise en m\'emoire quelconque.
Par exemple : 
\begin{itemize}
\item[] - Soient $u=\m\a[\a](f)\m\b[\a]x$ et $\th = \l x\l fu$. On a $val(u)=\{ 0 \}$ et
$rep(u)=\emptyset$.  \item[] - Soit $T=\l n((n)\l d\l g((T_i)n)\l x(g)( \so{s})(\so{p})x)\d$ o\`u \\
$\so{p}=\l n(((n)\l c((c)\so{0})\so{0})\l a\l c((c)(\so{s})(a)\l x\l yy)(a)\l x\l yx)\l x\l yy)\l x\l yy$ est
un $\l$-terme pour le "pr\'ed\'ecesseur" ($(\so{p})\so{0}\simeq\sb{\b}\so{0}$ et
$(\so{p})\so{n+1}\simeq\sb{\b}\so{n}$).
\end{itemize}
On v\'erifie que :
\begin{itemize}
\item[] - Pour tout $\th_n \simeq\sb{\b} \so{n}$, $((T)\th_n)f
\p \cases { (f)\so{0} &si $n=0$ \cr (f)(\so{s})(\so{p})(\so{s})^n\so{0}
&si $n\neq 0$ \cr}$
\\ Donc $T$ est un op\'erateur de mise en m\'emoire \footnote{Cet op\'erateur est
donn\'e par J.L. Krivine} pour pour les entiers.
\item[] -  $\v_{AF2} T:\q x\{N$*$[x] \f \neg\neg N[x] \}$ (le typage n\'ecessite l'introduction d'un
syst\`eme d'\'equations ad\'equat pour le type des entiers). 
\item[] - Mais $((T)\th)f \p
(f)(\so{s})(\so{p})\so{0}$ et $(\so{s})(\so{p})\so{0}\simeq\sb{\b}\so{1}$. $\Box$ 
\end{itemize}

\end{document}